\documentclass[a4paper,oneside,10pt]{article}

\usepackage[english]{babel}    
\usepackage[utf8]{inputenc}    
\usepackage[T1]{fontenc}       
\usepackage[reqno]{amsmath}    
\usepackage{amssymb,amsthm}    
\usepackage{mathtools}         
\usepackage{url}               
\usepackage{icomma}            
\usepackage{units}             
\usepackage{enumerate}         
\usepackage{array}             
\usepackage[usenames]{color}   
\usepackage{graphicx}          
\usepackage{caption}           
\usepackage{times}             
\usepackage{titlesec}          
\usepackage{soul}

\usepackage[bookmarks, colorlinks=true, linkcolor=black, anchorcolor=black,
citecolor=black, filecolor=black, menucolor=black, runcolor=black,
urlcolor=black, pdfencoding=unicode]{hyperref}  
\usepackage[
paper=a4paper,
top=2.5cm,
bottom=2.5cm,
left=2cm,
right=2cm
]{geometry}  

\renewcommand{\b}{\boldsymbol}
\newcommand{\p}{\b p}

\newcommand{\lap}{\nabla^2}

\newcommand{\dd}{\mathrm{d}}

\newcommand{\code}[1]{\texttt{#1}}

\newenvironment{itemize*}{\vspace{-6pt}\begin{itemize}\setlength{\itemsep}{0pt}\setlength{\parskip}{2pt}}{\end{itemize}}
\newenvironment{enumerate*}{\vspace{-6pt}\begin{enumerate}\setlength{\itemsep}{0pt}\setlength{\parskip}{2pt}}{\end{enumerate}}
\newenvironment{description*}{\vspace{-6pt}\begin{description}\setlength{\itemsep}{0pt}\setlength{\parskip}{2pt}}{\end{description}}

\newcommand{\Title}{GPU accelerated RBF-FD solution of Poisson's equation}
\newcommand{\Author}{Mitja Jančič, Jure Slak, Gregor Kosec}

\title{\Title}
\author{\Author}
\date{\today}
\hypersetup{pdftitle={\Title}, pdfauthor={\Author}, pdfcreator={\Author},
pdfproducer={\Author}, pdfsubject={}, pdfkeywords={}}  

\usepackage{multicol}
\usepackage{nth}

\newenvironment{Figure}
{\par\medskip\noindent\minipage{\linewidth}}
{\endminipage\par\medskip}
\newenvironment{Table}
{\par\medskip\noindent\minipage{\linewidth}}
{\endminipage\par\medskip}

\renewcommand{\b}{\boldsymbol}
\titleformat{\section}{\scshape\normalsize\centering}{\thesection.}{0.5em}{}
\titleformat{\subsection}{\itshape\normalsize}{\thesubsection.}{0.5em}{}
\titleformat{\subsubsection}{\itshape\normalsize\normalsize}{\thesubsection)}{0.5em}{}
\renewcommand{\thesection}{\Roman{section}}
\renewcommand{\thesubsection}{\Alph{subsection}}

\titlespacing\section{0pt}{2.82mm plus 2mm minus 0mm}{-0.71mm}
\titlespacing\subsection{0mm}{2.12mm}{-1.06mm}
\titlespacing\subsubsection{3.18mm}{0.72mm}{-2.12mm}
\titlespacing\paragraph{5.08mm}{0.0mm}{1ex}
\renewcommand{\thesection}{\Roman{section}}
\renewcommand{\thesubsection}{\Alph{subsection}}

\usepackage{soulutf8,color}

\usepackage{xcolor}
\usepackage{listings}

\colorlet{mygray}{black!30}
\colorlet{mygreen}{green!60!blue}
\colorlet{mymauve}{red!60!blue}

\begin{document}


\begin{center}
\fontsize{24}{28}\selectfont
\Title \\[2ex]
\fontsize{11}{11}\selectfont
Mitja Jančič$^{1,2}$, Jure Slak$^{1, 3}$, Gregor Kosec$^1$ 
\\[0.5mm]
\fontsize{10}{10}\selectfont
$^1$ ``Jožef Stefan'' Institute, Parallel and Distributed Systems 
Laboratory, Ljubljana,
Slovenia \\[0.5mm]
$^2$ ``Jožef Stefan'' International Postgraduate School, Ljubljana,
Slovenia \\[0.5mm]
$^3$ Faculty of Mathematics and Physics, University of Ljubljana, 
Ljubljana, Slovenia \\[1mm]
\href{mailto:mitja.jancic@ijs.si}{mitja.jancic@ijs.si},
\href{mailto:jure.slak@ijs.si}{jure.slak@ijs.si},
\href{mailto:gregor.kosec@ijs.si}{gregor.kosec@ijs.si}
\end{center}

\vspace{1ex}

\begin{multicols}{2}

\fontsize{9}{10}\selectfont
{\bfseries\noindent
	
	Abstract -- The Radial Basis Function-generated finite differences 
	became a
	popular variant of local meshless strong form methods due to its 
	robustness
	regarding the position of nodes and its controllable order of 
	accuracy.
	In this paper, we present a GPU accelerated numerical solution of 
	Poisson's 
	equation on scattered nodes in 2D for orders from 2 up to 6.  We 
	specifically 
	study
	the effect of using different orders on GPU acceleration efficiency.
}
\fontsize{10}{11}\selectfont

\section{Introduction}
\label{sec:intro}

In contrast to the traditional numerical methods for solving partial 
differential equations (PDE) that require connectivity between 
discretization 
nodes, meshless methods can operate on scattered nodes. Although 
seemingly 
minimal difference, this feature made meshless methods 
popular~\cite{fornberg2015solving} in various fields of science and 
engineering 
ranging from computational fluid dynamics~\cite{kosec2018local, 
	budiana2020meshless,zhang2018gpu} to option 
	pricing~\cite{rad2015pricing}. 
The key to becoming so popular is that discretization of arbitrary 
domain with 
scattered nodes is considered much easier problem in comparison with 
meshing. To some degree it can also be automated in dimensionless 
sense~\cite{slak2018generation}.  

From a historical point of view, meshless methods were introduced in 
1990s.
Since then, many meshless methods have been developed, e.g.\ the 
Element Free Galerkin
Method~\cite{Belytschko1994}, the Diffuse Element 
Method~\cite{Nayroles1992},
the Partition of Unity Method~\cite{MELENK1996289}, the Local 
Petrov-Galerkin
Methods~\cite{atluri1998new}, the h-p Cloud Methods~\cite{duarte1996h}, 
etc.

The radial basis function-generated finite differences (RBF-FD) was 
first
mentioned in~\cite{tolstykh2003using} as a local strong form meshless 
method
for solving PDEs. The method approximates differential operators only 
using
scattered nodes. However, often used RBFs, e.g.\ Gaussians, include a 
shape
parameter that can be crucial to the overall method
stability~\cite{flyer2016role}. Stability issue has been recently 
addressed by
using Polyharmonic splines (PHS) and additionally augmenting them with
polynomials~\cite{bayona2017role}.

The generality of the meshless methods comes with the price, namely 
higher 
complexity. First, the shape functions of stencil weights have to be 
computed every time from scratch. Second, the support sizes are 
typically much bigger 
than in mesh-based methods, especially in high order RBF-FD methods. 
The 
natural way to accelerate most computationally demanding 
parts of the 
solution procedure is to employ parallel computing. There have been 
several 
reports on parallel meshless solution procedures in the past, including 
distributed computing as well as shared memory 
approaches~\cite{zhang2018gpu, trobec2015parallel, kelly2014numerical, 
kosec2013local, 
	bollig2012solution}. In this paper we analyze the graphics 
	processing unit (GPU) 
acceleration of RBF-FD explicit solution of Poisson problem with 
Dirichlet
boundary condition. We are especially interested in the effect of
using different polynomial orders. 


%

The rest of the paper is organized as follows: in 
section~\ref{sec:meshless} a
short presentation of local meshless methods is given, in
section~\ref{sec:GPU} we present our workflow and how GPU was included 
in our
computations, in section~\ref{sec:problem} problem is explained, in
section~\ref{sec:results} results are presented, and in
section~\ref{sec:conclusions} final conclusions are given.

\section{Local strong form meshless methods}
\label{sec:meshless}

In general, the idea of meshless methods is the use of local 
discretization
points to construct an approximation of the considered field and later 
use it for
manipulation with differential operators. Discretization points are 
referred
to as computational nodes or just \emph{nodes}. The nodes are placed 
within the domain
and on its boundary, and their distribution can be scattered, as 
presented in
Fig.~\ref{fig:solution}.

\begin{Figure}
	\centering
	\includegraphics[width=0.99\linewidth]{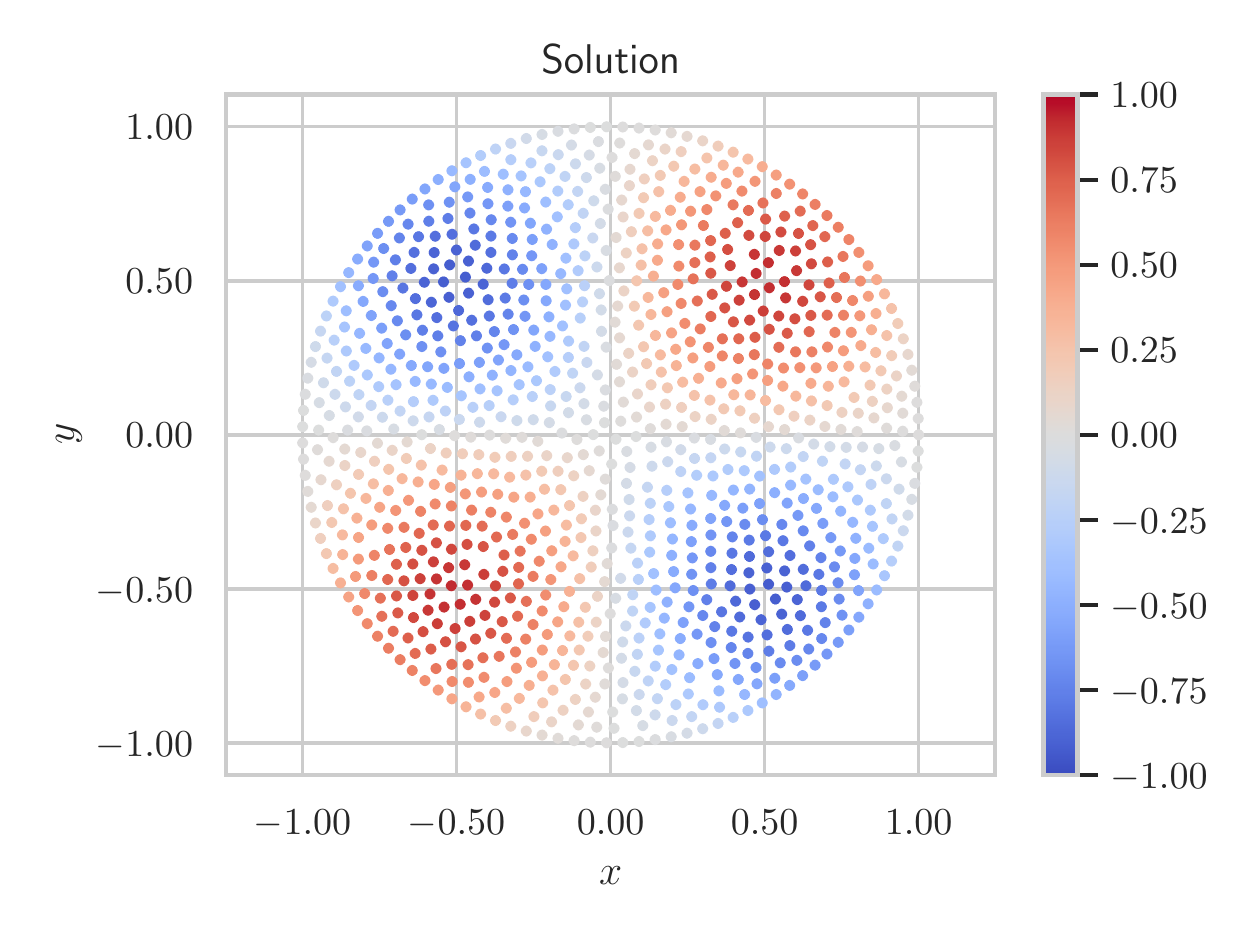}
	\captionof{figure}{Solution of Poisson's problem with Dirichlet 
	boundary
		conditions on scattered 2D nodes. Chosen highest polynomial 
		degree 
		$m=2$,
		support size $n = 15$ and number of nodes $N = 1027$.}
	\label{fig:solution}
\end{Figure}

Local strong form meshless methods approximate derivatives in the form
\begin{equation}
(\mathcal{L} u)(x) \approx \sum_{p \in N(x)} w_i^{\mathcal{L}}(x) u(p),
\end{equation}
where $N(x)$ is a set of neighboring nodes of $x$ and $\mathcal{L}$ is 
a differential operator.
Different ways of computing the weights $w$ exists. We will use the 
RBF-FD
method~\cite{slak2019high}.

We have implemented the RBF-FD based solution procedure using object
oriented approach and C++'s strong template system. Node positioning, 
support
selection, differential operator approximation and PDE discretization 
and other
modules are all available in one package, the Medusa 
library~\cite{medusa}.
Please refer to our open source Medusa library for more features and 
examples.

\section{GPU implementation}
\label{sec:GPU}
Compared to a CPU-based system, the memory bandwidth in GPU based systems is often an order 
of magnitude
higher, however at the price of higher memory latency. A single access to 
on board GPU
memory takes approximately 400-600 cycles compared to a floating point 
operation which takes
approximately 1-2 cycles. Current GPUs support 
large
numbers of processing elements and must be programmed in a high 
data-parallel
way in order to observe the advantages of GPU programming and to 
help
hide the memory latency. Large data is also needed to keep the 
processing
units busy and good knowledge of GPU systems is required to reduce the 
memory
communication and optimize the kernel functions~\cite{kosec2013local}.

In this work only part of the problem was dumped to the GPU. The 
illustration of
the GPU implementation is in Fig.~\ref{fig:scheme}. Green boxes are 
specific
to GPU programming while the yellow box presents the actual 
calculations done on
the GPU. Shapes, node positions and support nodes were still tasks 
executed by
the CPU, only the explicit time loop (yellow box in Fig.~\ref{fig:scheme}) was
ported to the GPU. Time loop execution was also timed using high resolution timer, enabling us to
estimate the speedups and thus evaluate the performance of CPU and GPU.

\begin{Figure}
	\centering
	\includegraphics[width=0.9\linewidth]{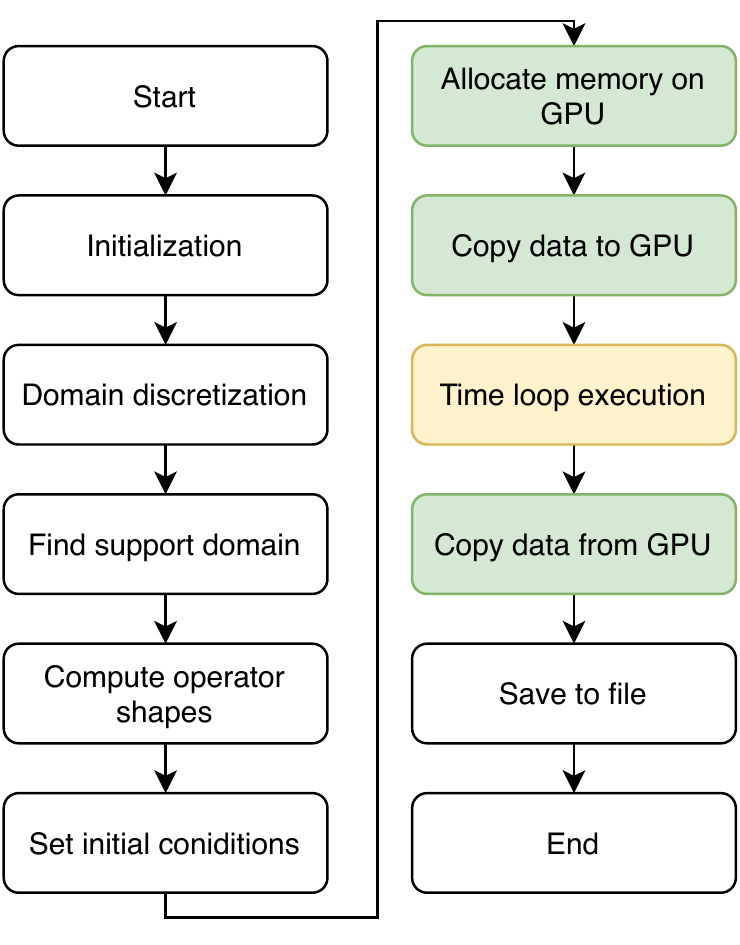}
	\captionof{figure}{The implementation scheme.}
	\label{fig:scheme}
\end{Figure}

\subsection{Time loop implementation}

Once the shapes, node positions and support nodes computation was carried out on the 
host (CPU), memory was allocated on the device (GPU) using 
\code{cudaMalloc} memory allocation function. As necessary data 
structures for numerically feasible problem sizes on the on-board GPU 
memory, we simply copied relevant data to the host using the 
\code{cudaMemcpy} function with \code{cudaMemcpyHostToDevice} argument 
and left the data there during time loop execution.

Time loop computation on the GPU was executed by introducing a 
\code{solveOnGpu} kernel function as shown in code 
listing~\ref{lst:code}.

\begin{lstlisting}[caption={Sample code of time loop execution on GPU.},captionpos=b, label={lst:code}]
	int tpb = 32;
	int bpg = (N + tpb - 1) / tpb;
	
	for (int step = 0; step < timeSteps; ++step) {
		solveOnGpu <<<bpg, tpb>>> (d_u1, d_u2, ...);
		cudaMemcpy(d_u1, d_u2, sizeof(double) * N, 
										cudaMemcpyDeviceToDevice);
	}
\end{lstlisting}

In Compute Unified Device Architecture (CUDA) the actual use of 
multi-processors 
and processing units is controlled by block 
size~\cite{nickolls2008scalable}. This reflects in the number of 
threads which act on each processor. In order to fully exploit the 
benefits of a GPU,
proper definition of block size is of great importance. When chosen to 
small,
processing units tend to stall, if chosen to large other effects (e.g.\ 
register
spilling) might interfere with top performance~\cite{kosec2013local}. 

Several tests of block size of multiples of 32 were made. Best time 
performance showed at $N_{tpb} = 32$ threads per block and $N_{bpg} = 
(N + N_{tpb}  - 1) / N_{tpb} $ blocks per grid with $N$ nodes in the 
domain. 

Once the time loop execution on the GPU is finished, \code{cudaDeviceSynchronize} function is 
used for synchronization and results are copied from the device to 
host using \code{cudaMemcpy} function with the 
\code{cudaMemcpyDeviceToHost} argument. The results are saved to a 
file and post-processing is done with python 3. To prevent memory 
leakage, \code{cudaFree(varibaleName)} was used before the shutdown.

\section{Problem description}
\label{sec:problem}
With the aim to analyze the performance of CPU and GPU in terms of 
support size
and highest augmented polynomial degree, proposed solution procedure 
and its
implementation is studied on a Poisson problem with Dirichlet boundary
conditions:
\begin{align}
\frac{\partial}{\partial t}u(\p) - \lap u(\p) &= f(\p) \quad &&\text{in 
} \Omega, \\
u(\p) &= \prod _{i = 1} ^d \sin(\pi _i) \quad && \text{on } \partial 
\Omega,
\end{align}
where $f(\p) = d \pi ^2 \prod_{i=1} ^d \sin(\pi p_i)$ and domain space 
$\Omega$ is a $d=2$ dimensional unit disk with boundary $\partial 
\Omega$.

The problem is discretized in both time and spatial domain space
\begin{equation}
\label{eq:explicit}
u_2(\p) = u_1(\p) + \dd t(f(\p) + \lap u_1(\p)),
\end{equation}
where $\dd t$ is an infinitesimal time step and $\lap $ is approximated 
as described in section~\ref{sec:meshless}.

The closed form solution of the above problem is $u(\p) = \prod_{i=1}^d 
\sin(\pi
p_i)$ allowing us to validate the numerically obtained solution.

\section{Results}
\label{sec:results}

Numerical results are computed using RBF-FD with PHS radial basis 
functions
$\Phi(r) = r^3$ and monomial augmentation. Radial function was kept the 
same for
all cases. The convergence order is directly controlled by the highest 
augmented
polynomial degree $m \in  \left \{ 2, 4, 6\right \} $, however the 
larger the
polynomial degree, the larger the recommended support size $n = 
\binom{m +
	d}{m}$. Recommended support size in terms of domain space 
	dimensionality and
highest polynomial degree was first given by 
Bayona~\cite{bayona2017role},
however larger supports can also be used~\cite{jani2019analysis}. In 
our test $n
\in  \left \{ 12, 15, 20, 30, 45, 60\right \} $.

All computations were performed on a single core of a computer with
\texttt{Intel(R) Xeon(R) CPU E5-2620 v3 @ 2.40GHz} processor and 64 GB 
of DDR4
memory and graphical processing unit \texttt{NVIDIA GeForce RTX 2080 
Ti, 11GB
	GDDR6} with \texttt{4352 CUDA cores}. Code was compiled using 
	\texttt{g++ (GCC)
	8.1.0} on Linux with \texttt{-O3 -DNDEBUG -std=c++11} flags while 
	GPU code was
compiled using \texttt{CUDA release V9.1.85}.

An example of numerical solution is shown in Fig.~\ref{fig:solution}.

\subsection{Execution times}
The discretized Poisson problem (\ref{eq:explicit}) was solved for 
$10^5$
time steps with time step $\dd t = 10^{-6}$ seconds, resulting in total 
simulation time 
$t = 0.1$ seconds. 

The polynomial degree takes an important role in the shape computation 
stage.
After shapes are stored the differential operator only
depends on the support. As illustrated on implementation scheme from
Fig.~\ref{fig:scheme}, the shapes and supports are in our case 
calculated on
the CPU. Different time loop execution times for various polynomial 
degrees $m$
are therefore not expected neither are they observed.

In Fig.~\ref{fig:times} we observe how the performance from CPU-only 
code is
directly proportional to number of nodes, while multiple regimes can be 
seen in
the single-GPU implementation.

\begin{Figure}
	\centering
	\includegraphics[width=0.99\linewidth]{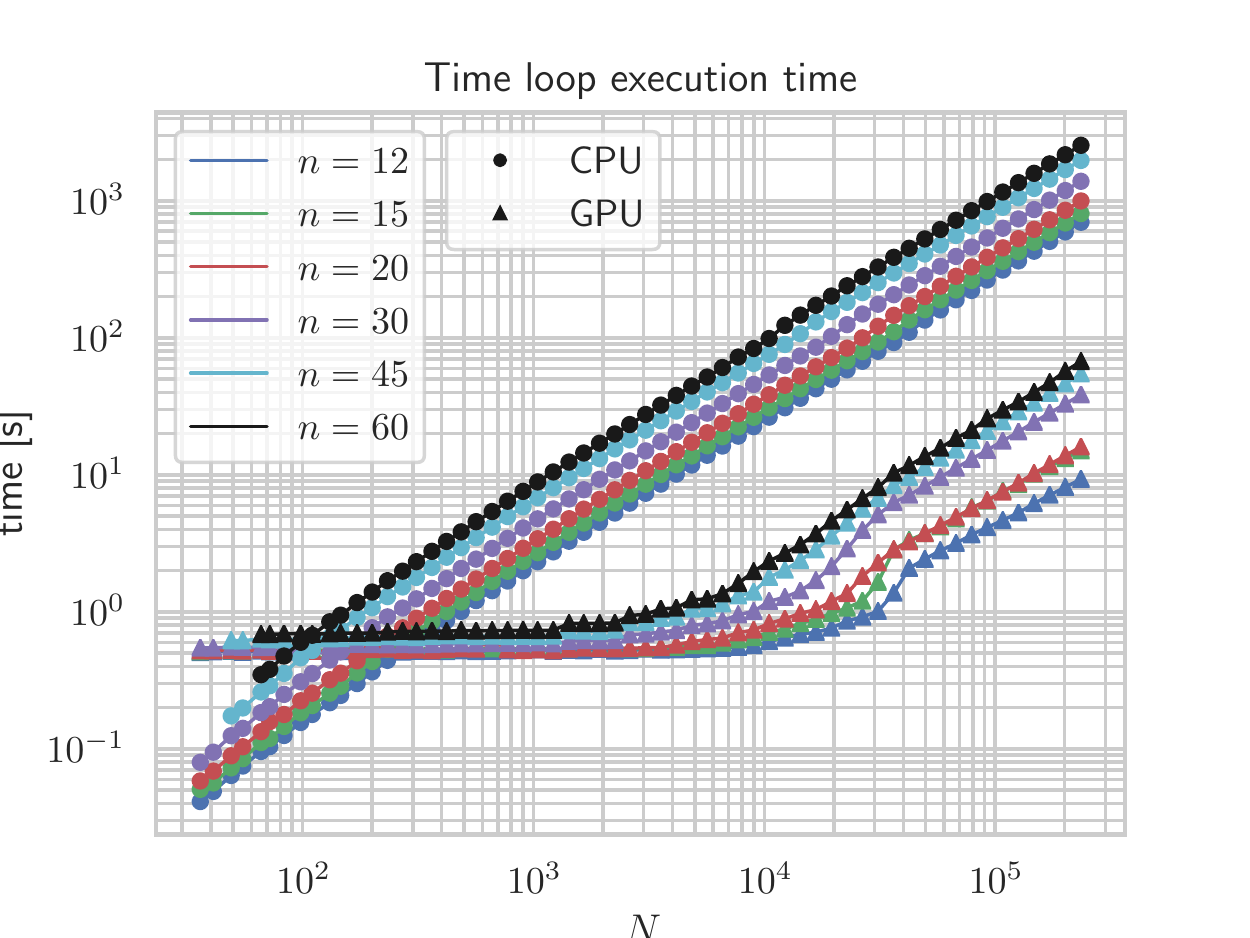}
	\captionof{figure}{Time loop execution times on single CPU and GPU. 
	Chosen 
		highest polynomial degree $m=2$.}
	\label{fig:times}
\end{Figure}

In the first regime it makes no significant difference what portion of 
the
available L2 cache is occupied by the data set. There is just not  
enough data
dumped to the GPU to fully exploit the advantages of parallelization as 
most of
the time is spent on communication and data distribution rather than on
computation itself. For larger data sets ($N>10^3$), the computation 
part prevails and
faster computation times are observed. Increasing data set size even 
further again leads to computation times proportional to the number of 
nodes - similar
to the CPU. As expected generally higher execution times are observed 
for larger
support sizes.

Speedup is defined as ratio between the execution of the time loop on a 
single
CPU $t_{CPU}$ and single GPU $t_{GPU}$. Speedups are shown on
Fig.~\ref{fig:speedup}. Similar regimes as in Fig.~\ref{fig:times} can 
be
seen starting with low or zero benefit from GPU implementation. 
Increasing the
data exploits the parallelization advantages which peak at certain data 
set size
dumped to the GPU.

\begin{Figure}
	\centering
	\includegraphics[width=0.99\linewidth]{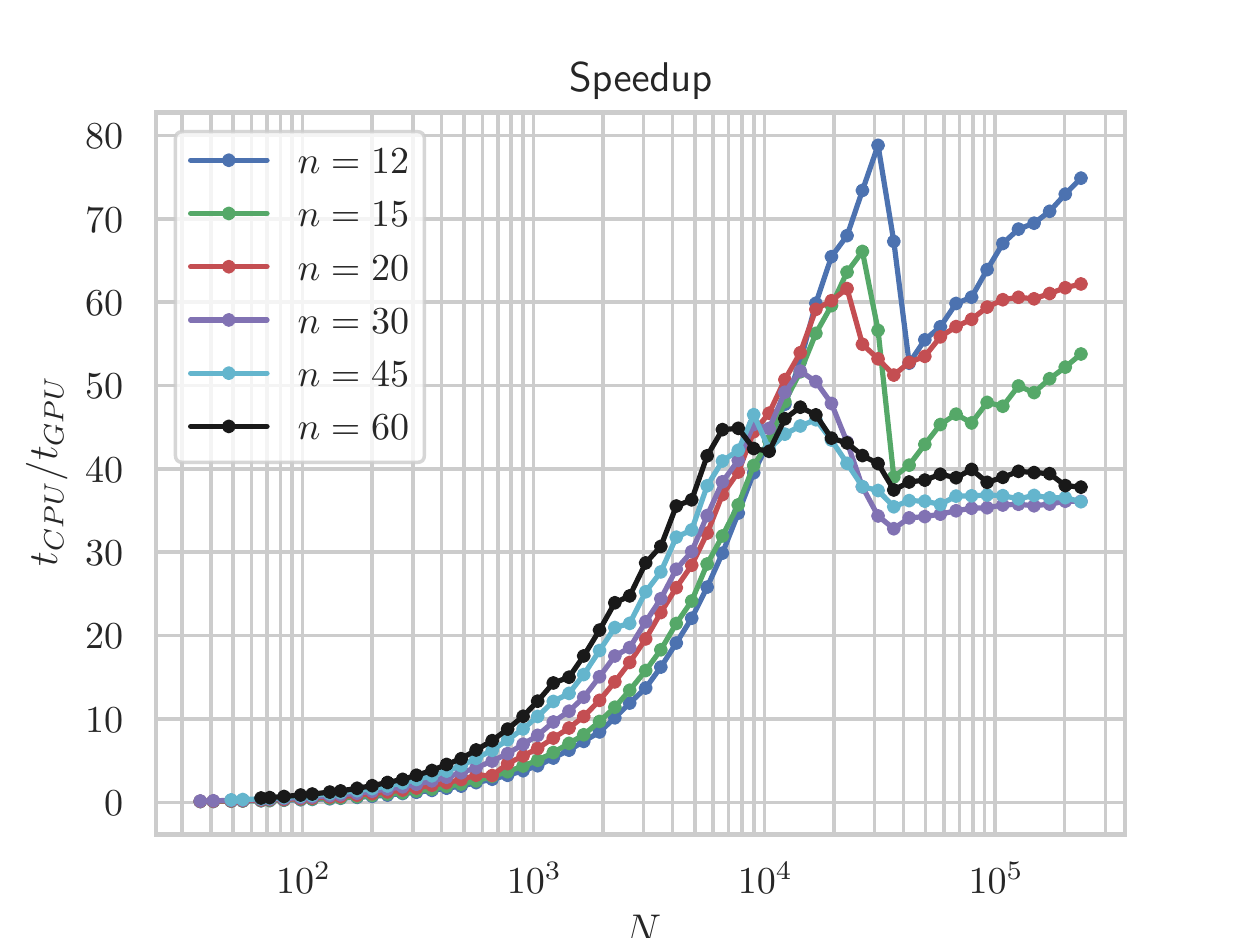}
	\captionof{figure}{Observed speedups.}
	\label{fig:speedup}
\end{Figure}

To explain the notable drop of speedup in Fig.~\ref{fig:speedup}, we 
have to
understand hardware structure in more detail. The L2 cache size of the
\texttt{NVIDIA GeForce RTX 2080 Ti} graphics processing unit is 5.5 
MB~\cite{RTX}.
Therefore increasing the data set size exploits the advantages of GPU
implementation but only until the L2 cache of the GPU is almost, but 
not yet,
full. Once the L2 cache is filled, more memory communication with higher
latency is needed.

A rough estimation of how much data is copied to the GPU can be made. 
In our
implementation we used doubles and integers, $N(4+n)$ double-precision 
and $N_i
+ Nn$ integer numbers, where $N_i$ is number of nodes in interior. 
Knowing that
size of double and integer is 8 and 4 bytes respectively, we can 
estimate when
the L2 cache is full and approximately estimate the speedup peak 
position. The
calculations are gathered in table~\ref{tab:esti}. Note that the 
estimated $N_{5.5
	\text{ MB}}$ is calculated by enforcing $N_i = N$. The larger the 
	number of
nodes, the more acceptable this enforcement is. Estimated data size 
copied on
GPU is also shown in Fig.~\ref{fig:est}.

\begin{Table}
	\centering
	\renewcommand{\arraystretch}{1.2}
	\begin{tabular}{rrr} \hline
		Support size & Estimated $N_{5.5\text{ MB}}$& $N$ speedup peak 
		\\
		\hline \rule{0pt}{12pt} 
		12 & 30556 & 31085 \\
		15 & 25463 & 26601 \\
		20 & 19928 & 19527 \\
		30 & 13889 & 14313 \\
		45 & 9549 & 8991 \\
		60 & 7275 & 7707 \\
		\hline
	\end{tabular}
	\captionof{table}{Estimated and observed peak performance as 
	function of number of nodes $N$.}
	\label{tab:esti}
\end{Table}

\begin{Figure}
	\centering
	\includegraphics[width=0.99\linewidth]{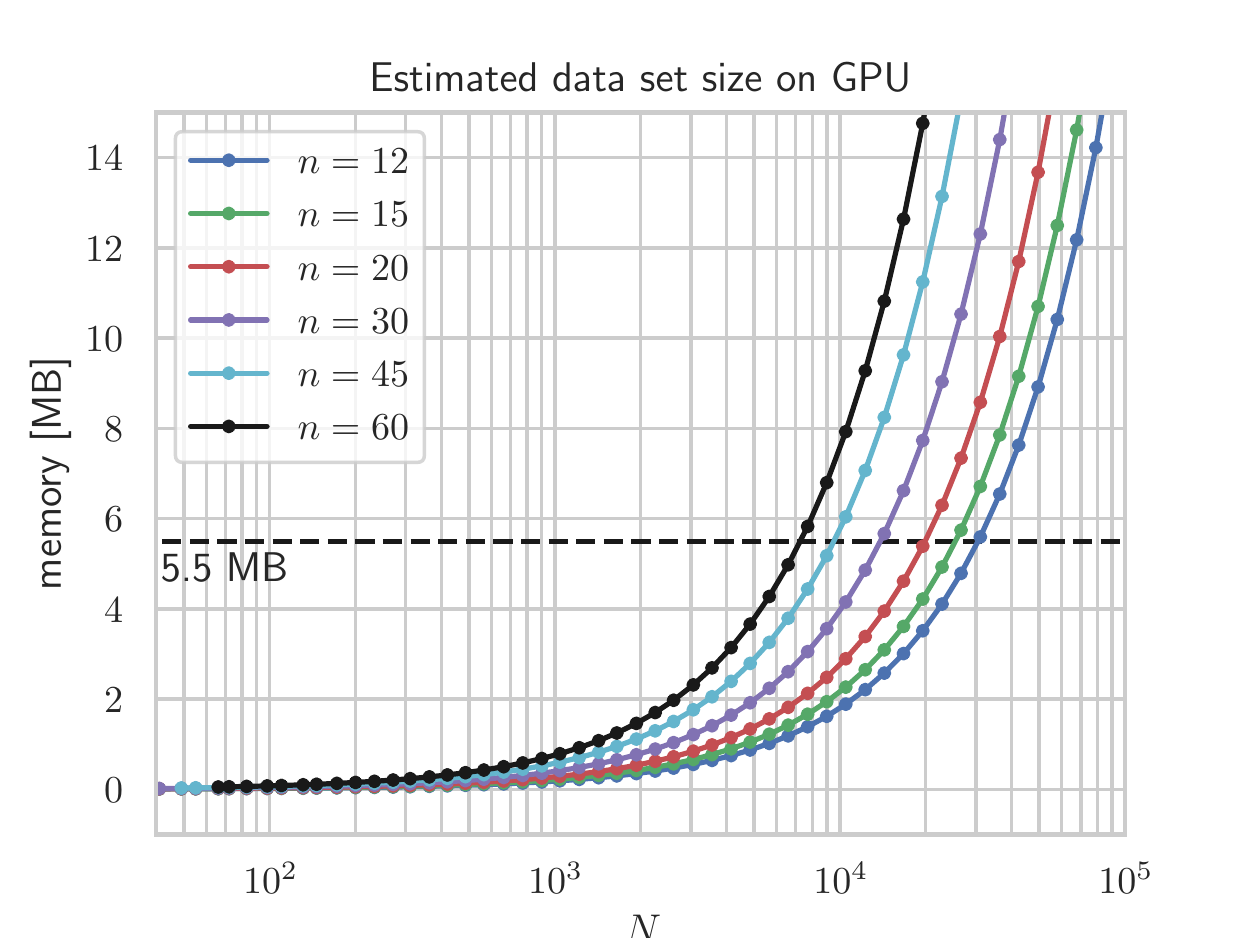}
	\captionof{figure}{Estimated data set size copied to GPU.}
	\label{fig:est}
\end{Figure}

Table~\ref{tab:esti} shows how observed speedup peaks are close to the 
estimated
calculations. We therefore conclude that the drop of speedup
performance is closely related to the size of L2 cache.

\section{Conclusions}
\label{sec:conclusions}

The execution performance of solution of a Poisson problem on 2D 
scattered 
nodes with Dirichlet boundary
conditions was analyzed in this paper. We measured how changing support 
sizes, 
a crucial parameter in using high order RBF-FD,
and total number of discretization nodes affects the parallel execution 
performance. We observed the speedup peaks, which can be explained by
measuring the data set size dumped to the GPU and comparing it with the
available L2 cache memory size.

In this paper only the execution of time loop was considered as 
performance
important, however some researchers already reported on using GPUs to 
calculate
the shapes. The next step could therefore be to dump even more 
computation to
the GPU. Additional extension could also be by solving the Poisson 
problem
implicitly or by comparing the computational performance using mulitple CPU and GPU devices.

\section*{Acknowledgments}
\label{sec:ack}
The authors would like to acknowledge the financial
support of the ARRS research core funding No.\ P2-0095
and the Young Researcher program PR-08346.

\bibliographystyle{unsrt}
\bibliography{ref}

\end{multicols}

\end{document}